\documentclass[12pt, a4paper]{article} 

\usepackage[margin=1in]{geometry} 

\usepackage{amsfonts}
\usepackage{amssymb}
\usepackage{amsmath}
\usepackage{amsthm}
\usepackage{dsfont}
\usepackage{graphics}
\usepackage{graphicx}
\usepackage{epsfig}
\usepackage{mathrsfs}
\usepackage[all]{xy}
\usepackage{color} 
\newtheorem{theorem}{Theorem}[section]
\newtheorem{assumption}[theorem]{Assumption}

\newtheorem{definition}[theorem]{Definition}

\newtheorem{lemma}[theorem]{Lemma}

\newtheorem{remark}[theorem]{Remark}

\numberwithin{equation}{section}

\usepackage{hyperref}
\hypersetup{
colorlinks=true, 
breaklinks=true, 
urlcolor= blue, 
linkcolor= red, 
citecolor=blue, 
pdftitle={}, 
pdfauthor={}, 
pdfsubject={} 
} 

\def\trieq{\triangleq}

\def\id{\mathbf{1}}
\def\qed{\hbox{ }\hfill {$\Box$}}
\def\ssd {\succeq_{(2)}}
\def\fsd {\succeq_{(1)}}
\def\isd{\succeq_{(i)}}

\def\proof{\noindent\textbf{Proof. }}

\def\VaR{\textrm{VaR}}

\def\ES{\textrm{ES}}

\DeclareMathOperator*{\mini}{minimize\,}

\def\Ebb{\mathbb{E}}

\def\Pbb{\mathbb{P}}

\def\Rbb{\mathbb{R}}

\def\Fc{\mathcal{F}}

\def\Rc{\mathcal{R}}

\def\Xc{\mathcal{X}}

\def\Asc{\mathscr{A}}

\def\Qsc{\mathscr{Q}}

\newcommand{\dd}{\ensuremath{\operatorname{d}\! }}

\newcommand{\ds}{\ensuremath{\operatorname{d}\! s}}
\newcommand{\dr}{\ensuremath{\operatorname{d}\! r}}
\renewcommand{\geq}{\geqslant}
\renewcommand{\leq}{\leqslant}
\renewcommand{\ge}{\geqslant}
\renewcommand{\le}{\leqslant}

\title{\large\textbf{Minimal Quantile Functions Subject to \\ Stochastic Dominance Constraints}\thanks{Supported by National Key R\&D Program of China (NO. 2018YFA0703900), NSFC (NO. 12071146, No.11971409, and No.12171169), Guangdong Basic and Applied Basic Research Foundation (Grant No. 2019A1515011338), The Hong Kong RGC (GRF No.15202421), The PolyU-SDU Joint Research Center on Financial Mathematics, The CAS AMSS-PolyU Joint Laboratory of Applied Mathematics, and The Hong Kong Polytechnic University.
}}

\author{Xiangyu Wang\thanks{School of Mathematical Sciences, University of Chinese Academy of Sciences, Beijing 100049, China; Email: \url{xywang93@hotmail.com}.}
\and
Jianming Xia\thanks{RCSDS, NCMIS, Academy of Mathematics and Systems Science, Chinese Academy
of Sciences, Beijing 100190, China; Email: \url{xia@amss.ac.cn}.}
\and
Zuo Quan Xu\thanks{Department of Applied Mathematics, The Hong Kong Polytechnic University, Kowloon, Hong Kong, China; Email: \url{maxu@polyu.edu.hk}.}
\and
Zhou Yang\thanks{School of Mathematical Sciences, South China Normal University, Guangzhou 510631, China; Email: \url{yangzhou@scnu.edu.cn}.}
}

\begin{document}

\maketitle

\begin{abstract}
We consider a problem of finding an SSD (second-order stochastic dominance)-minimal quantile function subject to the mixture of FSD (first-order stochastic dominance) and SSD constraints. The SSD-minimal solution is explicitly worked out and has a close relation to the Skorokhod problem. This result is then applied to explicitly solve a risk minimizing problem in financial economics. 
\medskip


\end{abstract}

\section{Introduction} 
Stochastic dominance (first/second-order) plays an important role in statistics, financial economics, insurance, etc. In this paper we consider a problem of finding an SSD (second-order stochastic dominance)-minimal quantile function subject to the mixture of FSD (first-order stochastic dominance) and SSD constraints. 
More precisely, given two quantile functions $Q_1$ and $Q_2$, the problem is to find a quantile function $Q^*$ such that:
\begin{itemize}
\item[(i)] $Q^*$ first-order stochastic dominates $Q_1$ and second-order stochastic dominates $Q_2$;
\item[(ii)] $Q^*$ is the SSD-minimal one, in the sense that if a quantile function $Q$ also satisfies condition (i), then $Q$ second-order stochastic dominates $Q^*$.
\end{itemize} 
{ Such an SSD-minimal solution has applications to determine the lower optimal stopping value of a sequence of random variables whose joint distribution is not known (only the marginal distributions are known); see M\"uller and R\"uschendorf \cite{MR01}.}
For the special case when the benchmark distribution functions have 
at most countable crossing points, M\"uller and R\"uschendorf \cite{MR01} provided an explicit construction of the SSD-minimal solution. In the iterated application of the SSD-minimal solution of M\"uller and R\"uschendorf \cite[Proposition 2.2 and Theorem 2.1]{MR01}, it is, however, not easy to verify that $F_i$ and $G_{i+1}^*$ have at most countable crossing points since $G_{i+1}^*$ itself is not a given distribution but a part of the solution of the problem. 
{To facilitate the applications of SSD-minimal quantile functions, we need to find an explicit construction for the general case.} Such an observation motivates us to discuss the SSD-minimal problem for the general case.

We will study the general case and show that the problem is highly related to the Skorokhod problem.
The main contribution of this paper is to show such an SSD-minimal quantile function always exists in the general case and to provide an explicit expression for the unique solution. 

Apart from the problem of lower optimal stopping value of random variables in M\"uller and R\"uschendorf \cite{MR01}, the SSD-minimal quantile function can also be applied to explicitly solve a risk minimizing problem in financial economics, which reduces to an EMP (expenditure minimizing problem). 

The classical EMP is to find an optimal payoff which has a given probability distribution function and minimizes the price of its replication. This problem can go back at least to Dybvig \cite{D88} and is further investigated by Schied \cite{S04}, Carlier and Dana \cite{CD06}, and Jin and Zhou \cite{JZ08}. The related result plays an important role in the so called ``quantile formulation" of the problem to find the optimal payoffs for preferences described by non-expected utilities such as rank-dependent utilities, rank linear utilities, and cumulative prospect theory; see, e.g.,
Carlier and Dana \cite{CD06}, Jin and Zhou \cite{JZ08}, He and Zhou \cite{HZ11},
Xia and Zhou \cite{XZ16}, and Xu \cite{X16}.

In a non-atomic probability space and a complete market, the EMP with an FSD constraint is essentially equivalent to the problem with the constraint that the payoff has a given distribution. 
Moreover, the problem with an SSD constraint
was investigated by F\"ollmer and Schied \cite{FS11} and has the same solution as the one with an FSD constraint; see 
F\"ollmer and Schied \cite[Theorem 3.44]{FS11}. 
The risk minimizing problem in this paper is related to the EMP with a mixture of an FSD constraint and an SSD constraint. We will investigate this problem and provide an explicit optimal solution.

The rest part of the paper is organized as follows. Section \ref{sec:main} formulates a problem of finding an SSD-minimal quantile function subject to the mixture of FSD and SSD constraints and presents our main result. 
Section \ref{sec:emp} provides an application of the main result to a risk minimizing problem in financial economics. The main result is proved in Section \ref{sec:proof}.

\section{Problem Formulation and Main Result}\label{sec:main}

We will use quantile functions,\footnote{Quantile functions are the general inverse of probability distribution functions.} instead of probability distribution functions, to formulate the problem. 
Let $\Qsc$ denote the set of (upper) quantile functions of all probability distribution functions, 
then
$$\Qsc=\big\{Q: (0,1)\to\Rbb \,\big|\, Q \mbox{ is increasing and right-continuous}\big\}.\footnote{Throughout the paper ``increasing" means ``non-decreasing" and ``decreasing" means ``non-increasing." }$$
For more details about quantile functions,
see F\"ollmer and Schied \cite[Appendix A.3]{FS11}.

Let $Q$, $Q_0\in\Qsc$. We say that $Q$ first-order stochastic dominates $Q_0$ and write it $Q\fsd Q_0$, if 
$$Q(t)\ge Q_0(t) \quad\mbox{for all }t\in(0,1).$$
We say that $Q$ second-order stochastic dominates $Q_0$ and write it $Q\ssd Q_0$, if
$$\int_0^tQ(s)\ds \geq \int_0^tQ_0(s)\ds\quad\mbox{for all }t\in(0,1).$$
Clearly, $Q\fsd Q_0$ implies $Q\ssd Q_0$. 
Given a quantile function $Q_0\in \Qsc$, let 
$$\Qsc_i(Q_0)\trieq\left\{Q\in\Qsc\,\left|\,Q\isd Q_0\right.\right\},\quad i\in\{1,2\}.$$
Both $\Qsc_1(Q_0)$ and $\Qsc_2(Q_0)$ are obviously convex sets. 

{
Consider two benchmark quantile functions $Q_1$, $Q_2\in\Qsc$ that satisfy 
\begin{assumption}\label{ineq:q12-}
$\int_0^1Q_i^-(s)\ds<\infty,\quad i\in\{1,2\}.$
\end{assumption}
}

\begin{definition}
A quantile function $Q^*\in\Qsc_1(Q_1)\cap\Qsc_2(Q_2)$ is called \textit{SSD-minimal} in
$\Qsc_1(Q_1)\cap\Qsc_2(Q_2)$ if $Q\in\Qsc_2(Q^*)$ for every $Q\in\Qsc_1(Q_1)\cap\Qsc_2(Q_2)$. 
\end{definition}

One can easily show that there is at most one SSD-minimal quantile function in
$\Qsc_1(Q_1)\cap\Qsc_2(Q_2)$. 
M\"uller and R\"uschendorf \cite{MR01} gave an explicit construction of the SSD-minimal quantile function when the two benchmark distributions have at most countable discrete set of crossing points. 
We will show that, in the general case, the SSD-minimal quantile function $Q^*$ always exists and can be characterized
by the solution to the following ODE (ordinary differential equation)\footnote{{A solution to \eqref{opt:equa} is an absolutely continuous function $\phi:[0,1)\to\Rbb$ that satisfies \eqref{opt:equa}. }} with respect to $\phi$: 
\begin{align}\label{opt:equa}
\begin{cases}
\min\{\phi',\phi-f\}=0\quad\mbox{a.e. (almost everywhere) in } (0,1),\\
\phi(0)=0,
\end{cases}
\end{align}
where the function $f$ is defined as 
$$f(t)\trieq \int^t_0 \big(Q_2(r)-Q_{1}(r)\big)\dr \quad\mbox{for all }t\in [0,1).$$
{Obviously, Assumption \ref{ineq:q12-} guarantees that $f$ is finite and continuous on $[0,1)$.}

The next theorem characterizes the solution of ODE \eqref{opt:equa}, which, as we will see, is closely related to the Skorokhod problem.

\begin{theorem}\label{thm:ODEsolution}
{Under Assumption \ref{ineq:q12-},} the unique solution to ODE \eqref{opt:equa} is given by 
\begin{align}\label{defphi}
\phi(t)= \max_{0\leq s\leqslant t}f(s)\quad\mbox{for all } t\in [0,1).
\end{align}
\end{theorem} 

\proof 
ODE \eqref{opt:equa} is essentially the well-known Skorokhod problem; see, e.g., Revuz and Yor \cite[Lemma VI.2.1]{RY99}.
Actually, let $z(t)=\phi(t)-f(t)$, then \eqref{opt:equa} is obviously equivalent to 
\begin{description}
\item[(i)] $z\ge0$,
\item[(ii)] $\phi(0)=0$, $\phi$ is increasing, and
\item[(iii)] $\phi$ is flat off $\{t\in(0,1)\,|\, z(t)=0\}$, i.e., $\int_0^1\id_{z(t)>0}\dd\phi(t)=0$.\footnote{Hereafter, $\id_{A}$ is the indicator function for a statement $A$, whose value is 1 if the statement $A$ is true and $0$ otherwise.}
\end{description}
Therefore, by the Skorokhod lemma, the unique solution to ODE \eqref{opt:equa} is \eqref{defphi}. \qed

The next theorem characterizes the SSD-minimal quantile function $Q^*$,
whose proof is deferred to Section \ref{sec:proof}.

\begin{theorem} \label{thm:opt:eqtn}
{Suppose Assumption \ref{ineq:q12-} hold.}
Let $\phi$ be given in Theorem \ref{thm:ODEsolution}, and let
\begin{align}\label{defQ}
Q^{*}(t)\trieq Q_{1}(t)\id_{\phi(t)>f(t)}+(Q_{1}(t)\vee Q_{2}(t))\id_{\phi(t)=f(t)}\quad\mbox{for all } t\in(0,1). 
\end{align}
Then $Q^*$ is SSD-minimal in $\Qsc_1(Q_1)\cap\Qsc_2(Q_2)$.
\end{theorem} 

\begin{remark}
We know that $Q\in\Qsc_1(Q_1)\cap\Qsc_2(Q_2)$ if and only if
\begin{equation}\label{ineq:min:QQ12}
\min\left\{Q(t)-Q_1(t), \int_0^tQ(r)\dr-\int_0^tQ_2(r)\dr\right\}\ge0\quad\mbox{for all } t\in(0,1).
\end{equation}
Let $Q^*$ be given as in Theorem \ref{thm:opt:eqtn}. Then
by \eqref{opt:equa} we get 
$$\min\left\{Q^*(t)-Q_1(t), \int_0^tQ^*(r)\dr-\int_0^tQ_2(r)\dr\right\}=0\quad\mbox{for all } t\in (0,1).$$
Therefore, the SSD-minimal quantile function $Q^*$ makes the constraint in \eqref{ineq:min:QQ12} binding for all $t\in(0,1)$.
\end{remark}

\begin{remark}\label{mult} In the previous discussion, the problem has been investigated for the mixture of only one FSD constraint and only one SSD constraint. Actually, the multiple FSD (resp. SSD) constraints can be reduced to a single FSD (resp. SSD) constraint; see, e.g., Wang and Xia \cite[Remarks 3.5]{WX21} for details. 
\end{remark}

\section{Risk Minimizing Problem}\label{sec:emp}

In this section we apply our main result to a risk minimizing problem in financial economics.
Consider a non-atomic complete probability space $(\Omega, \Fc, \Pbb)$. 
Let $L^0$ denote all $\Fc$-measurable random variables {and $L^1$ all integrable $X\in L^0$}. 
For every $X\in L^0$, 
its (upper) quantile function is defined by 
\[Q_X(t)= \inf\{x\in\Rbb: \Pbb(X\leq x)>t\},\quad t\in (0,1).\]

Given two benchmark payoffs $X_1\in L^0$ {and $X_2\in L^1$}, consider the acceptance sets
$$\Asc_i=\big\{X\in L^0\,\big|\, X\isd X_i\big\},\quad i\in\{1,2\},$$
where $X\isd X_i$ is synonymous to $Q_{X}\isd Q_{X_i}$. 
Let the risk measures $\Rc_i$, $ i\in\{1,2\}$, be derived from $\Asc_i$ as follows:
$$\Rc_i(X)\trieq\inf\big\{m\in\Rbb\,\big|\, X+m\in\Asc_i\big\}{\mbox{ with }\inf\emptyset=\infty}, \quad X\in L^0.$$
Then $\Rc_1$ is called the \emph{loss value at risk} and $\Rc_2$ the \emph{benchmark adjusted expected shortfall}; see Bignozzi et al \cite{BBM20} and Burzoni et al \cite{BMW22}.
Moreover, we have the representations 
$$\Rc_1(X)=\sup_{t\in(0,1)}\{\VaR_t(X)-\VaR_t(X_1)\} $$
and 
$$\Rc_2(X)=\sup_{t\in(0,1)}\{\ES_t(X)-\ES_t(X_2)\},$$
where $\VaR_t(X)$ and $\ES_t(X)$ are the \textit{value at risk} and \textit{expected shortfall} of $X$ at the level $t\in(0,1)$, defined as
$$\VaR_t(X)\trieq -Q_X(t) \quad\mbox{ and}\quad \ES_t(X)\trieq -{1\over t}\int_0^tQ_X(s)\ds.$$ 
For a general relation between acceptance sets and monetary risk measures, see, e.g., F\"ollmer and Schied \cite{FS11}. 

It is natural to consider an acceptance set $\Asc$ given by
$$\Asc \trieq\big\{X\in L^0\,\big|\, X\isd X_i\mbox{ for }i\in\{1,2\}\big\}=\Asc_1\cap\Asc_2.$$
Let the risk measure $\Rc$ be derived from $\Asc$, that is, 
\begin{align}\label{riskmeasure}
\Rc(X)\trieq\inf\big\{m\in\Rbb\,\big|\, X+m\in\Asc_1\cap\Asc_2\big\}{ \mbox{ with }\inf\emptyset=\infty},\quad X\in L^0.
\end{align}
Clearly, {$\Rc$ is monotone and translation invariant,} $\Rc=\Rc_1\vee\Rc_2$, 
and $X\in\Asc_1\cap\Asc_2$ iff $\Rc(X)\le0$. 
{We proceed with the following technical assumption.
\begin{assumption}\label{ass:Q12rho}
$\int_0^1|Q_{X_2}(s)|\ds<\infty$ and $\int_0^1|Q_{X_i}(s)|Q_\rho(1-s)\ds<\infty$ for $i\in\{1,2\}$.
\end{assumption}}
Consider an arbitrage-free market. Assume the market is complete and has a
unique SDF (stochastic discount factor\footnote{Also sometimes termed
``pricing kernel" or ``state price density" in financial economics literature.}) $\rho\in L^0$ satisfying $ \Pbb(\rho>0)=1$
and $\Ebb[\rho]<\infty$. 
Given a budget level $x>0$, the risk minimizing problem for risk measure $\Rc$ is defined as 
\begin{align}\label{staticproblem}
\mini_{{X\in\Xc_\rho}} &\;\Rc(X)\;\; \mbox{ subject to }\;\Ebb[\rho X]\le x,
\end{align}
{where}
$${\Xc_\rho\trieq\{X\in L^0\,|\, \Ebb[\rho X^+]<\infty\}.}$$

Thanks to the monotonicity of the risk measure $\Rc$, solving problem \eqref{staticproblem} can reduce to solving the following EMP 
\begin{align}\label{expenditureminimizing}
\mini_{X\in\Xc_\rho} &\;\Ebb[\rho X]\;\;\mbox{ subject to }\Rc(X)\le 0.
\end{align} 
Lemma \ref{equivalent} below states quantitatively the relationship between the optimal solution to \eqref{staticproblem} and that to \eqref{expenditureminimizing}.

\begin{lemma}\label{equivalent} {Suppose Assumption \ref{ass:Q12rho} hold.
Then the optimal values of problems \eqref{staticproblem} and \eqref{expenditureminimizing} are finite. } If $X^*$ is an optimal solution to problem \eqref{staticproblem}, then
$X^*+\Rc(X^{*})$ is an optimal solution to problem \eqref{expenditureminimizing}. 
On the other hand, if $X^*$ is an optimal solution to problem \eqref{expenditureminimizing}, 
then $X^*+\frac{x-\Ebb[\rho X^*]}{\Ebb[\rho]}$ is an optimal solution to problem \eqref{staticproblem}. 
\end{lemma}

\proof
{Without loss of generality, we can assume that $X_1$ and $X_2$ are anti-comonotonic with $\rho$ because $\Rc$ is law-invariant. In this case, Assumption \ref{ass:Q12rho} reads $\Ebb[|X_2|]<\infty$ and $\Ebb[\rho|X_i|]<\infty$, $i\in\{1,2\}$.
Let $v_1$ and $v_2$ denote, respectively, the optimal values of problems \eqref{staticproblem} and \eqref{expenditureminimizing}.
Let $c$ be any real constant. Then 
\[ \quad \Rc_1(|X_1|+c)=\Rc_1(|X_1|)-c\le \Rc_1(X_1)-c=-c<\infty.\]
By the well known fact that $\ES_t(X_2)\ge -\Ebb[X_2]$, 
\begin{align*}
\Rc_2(|X_1|+c) &\le \Rc_2(c)=\sup_{t\in(0,1)}\{-c-\ES_t(X_2)\}\le -c+\Ebb[X_2]<\infty.
\end{align*}

Since $\Ebb[\rho]>0$, for sufficiently small $c$, we have that $\Ebb[\rho (|X_1|+c)]\leq x$ and $|X_1|+c$ is a feasible solution to problem \eqref{staticproblem}, which leads to $v_1\leq\Rc(|X_1|+c)<\infty$. 
Similarly, for sufficiently large $c$, we have that $\Rc(|X_1|+c)\leq 0$ and $|X_1|+c$ is a feasible solution to problem \eqref{expenditureminimizing}, which yields $v_2\leq \Ebb[\rho(|X_1|+c)]<\infty$.

On the other hand, 
using that $X_1$ is anti-comonotonic with $\rho$ and the Hardy--Littlewood inequality,
we have, for any $X\in \Xc_\rho$, 
\begin{align}
\Rc_1(X)\Ebb[\rho] &=\sup_{t\in(0,1)}(Q_{X_1}(t)-Q_{X}(t))\int_0^1 Q_{\rho}(1-s)\ds\nonumber\\
&\geq \int_0^1 (Q_{X_1}(s)-Q_{X}(s))Q_{\rho}(1-s)\ds \geq \Ebb[\rho X_1]-\Ebb[\rho X].\label{lower1}
\end{align}
For any feasible solution $X$ to problem \eqref{staticproblem}, it follows 
\[\Rc(X)\ge\Rc_1(X)\geq\frac{\Ebb[\rho X_1]-\Ebb[\rho X]}{\Ebb[\rho]}\geq\frac{\Ebb[\rho X_1]-x}{\Ebb[\rho]}.\] This together with $v_1<\infty$ implies that $v_1$ is finite.
For any feasible solution $X$ to problem \eqref{expenditureminimizing}, 
since $\Rc_1(X)\leq \Rc(X)\leq 0$, we see from \eqref{lower1} that 
\[\Ebb[\rho X] \geq \Ebb[\rho X_1].\] 
This together with $v_2<\infty$ implies that $v_2$ is finite.
}

Now suppose $X^*$ is an optimal solution to problem \eqref{staticproblem}. 
Then $d=\Rc(X^{*})$ is finite. 
Let $Y=X^*+d$. Then \[\Rc(Y)=\Rc(X^{*})-d=0.\] Suppose on the contrary that $Y$ is not optimal to \eqref{expenditureminimizing}, then there exists a $Z\in \Xc_\rho$ such that $\Ebb[\rho Z]< \Ebb[\rho Y]$ and $\Rc(Z)\leq 0$. In this case, as $\Ebb[\rho]<\infty$,
\[\Ebb[\rho (Z-d)]<\Ebb[\rho (Y-d)]=\Ebb[\rho X^*]\leq x.\] 
Hence, there exists a constant $\delta>0$ such that $\Ebb[\rho (Z-d+\delta)]\leq x$. Obviously, $Z-d+\delta\in\Xc_\rho$.
But 
\[\Rc(Z-d+\delta)=\Rc(Z)+d-\delta< d=\Rc(X^*),\]
which contradicts the optimality of $X^{*}$ to \eqref{staticproblem}. Hence, $X^*+d$ is an optimal solution to \eqref{expenditureminimizing}.

Suppose $X^*$ is an optimal solution to problem \eqref{expenditureminimizing}. 
Then $d=\frac{x-\Ebb[\rho X^*]}{\Ebb[\rho]}$ is finite. Let $Y=X^*+d$. Clearly 
\[\Ebb[\rho Y]=\Ebb[\rho X^*]+\Ebb[\rho d]=x.\]
Suppose on the contrary that $Y$ is not optimal to \eqref{staticproblem}, then there exists a $Z\in\Xc_\rho$ such that $\Ebb[\rho Z]\leq x$ and $\Rc(Z)<\Rc(Y)$. In this case, 
\[\Rc(Z-d)<\Rc(Y-d)=\Rc(X^*)\leq 0.\]
Hence, there exists a small constant $\delta>0$ such that \[\Rc(Z-d-\delta)=\Rc(Z-d)+\delta<0.\] Obviously, $Z-d+\delta\in\Xc_\rho$.
But 
\[\Ebb[\rho (Z-d-\delta)]=\Ebb[\rho Z]-d\Ebb[\rho]-\delta \Ebb[\rho]\leq x-(x-\Ebb[\rho X^*])-\delta \Ebb[\rho] <\Ebb[\rho X^*],\]
which clearly contradicts the optimality of $X^{*}$ to \eqref{expenditureminimizing}. Therefore, $X^*+d$ is an optimal solution to \eqref{staticproblem}.
\qed

By Lemma \ref{equivalent}, solving problem \eqref{staticproblem} is essentially equivalent to solving problem \eqref{expenditureminimizing}, so we focus on the latter from now on. 

As $\Rc$ is law-invariant and monotone, 
using the Hardy--Littlewood inequality, we know that any
solution $X^*$ to problem \eqref{expenditureminimizing} satisfies 
$X^*=Q_{X^*}(1-\xi)$ and the minimum is 
$$\Ebb[\rho X^*]=\int_0^1Q_{X^*}(s)Q_{\rho}(1-s)\ds,$$
where $\xi\in\Xi$ and\footnote{For the existence of such a $\xi$, see, e.g., F\"{o}llmer and Schied \cite[Lemma A.28]{FS11}. Moreover, let $\xi$ be a random variable uniformly distributed on $(0,1)$. Then by Xu \cite[Theorem 5]{X14}, $\xi\in\Xi$ if and only if $(\xi,\rho)$ is comonotonic.}
$$\Xi\trieq\{\xi\mid \xi\mbox{ is uniformly distributed on $(0,1)$ and }\rho=Q_\rho(\xi)\mbox{ a.s.}\};$$
see, e.g., Dybvig \cite{D88}, Schied \cite{S04}, Carlier and Dana \cite{CD06}, Jin and Zhou \cite{JZ08}, Xu \cite{X14}, and {Liebrich and Munari \cite{LM22}}. Recalling that $\Rc(X)\le0$ if and only if $X\in\Asc_1\cap\Asc_2$, we know that problem \eqref{expenditureminimizing} 
reduces to the following problem: 
\begin{align*} 
\mini_{X\in\Xc_\rho}\; \int_0^1Q_X(s)Q_\rho(1-s)\ds 
\quad\mbox{subject to } X\in\Asc_1\cap\Asc_2.
\end{align*}
Its quantile formulation is then given by 
\begin{align}\label{opt:Q:fsd+ssd}
\mini_{Q\in\Qsc_\rho}\; \int_0^1Q(s)Q_\rho(1-s)\ds 
\quad\mbox{subject to } Q\in\Qsc_1(Q_{X_1})\cap\Qsc_2(Q_{X_2}),
\end{align}
where
$$\Qsc_\rho\trieq\left\{Q\in\Qsc\,\left|\,\int_0^1Q^+(s)Q_\rho(1-s)\ds<\infty\right.\right\}.$$

By F\"ollmer and Schied \cite[Theorem 2.57 and Lemma 3.45]{FS11}, we have, for any $\overline{Q}\in\Qsc$,
\begin{equation}\label{eq:Q:Q0}
Q\in \Qsc_2(\overline{Q})\;\Longrightarrow\; \int_0^1Q(s)Q_\rho(1-s)\ds\ge \int_0^1\overline{Q}(s)Q_\rho(1-s)\ds.
\end{equation}
{Under Assumption \ref{ass:Q12rho},}\footnote{{Obviously, $\int_0^1|Q_{X_i}(s)|Q_\rho(1-s)\ds<\infty$ implies that $\int_0^1Q_{X_i}^-(s)\ds<\infty$, so $Q_i=Q_{X_i}$, $i\in\{1,2\}$, fulfill Assumption \ref{ineq:q12-}.}} let $Q^*$ be the SSD-minimal quantile function in $\Qsc_1(Q_{X_1})\cap\Qsc_2(Q_{X_2})$, which is given by Theorem \ref{thm:opt:eqtn}.
Then Assumption \ref{ass:Q12rho} and \eqref{eq:Q:Q0} imply that $Q^*$ solves problem \eqref{opt:Q:fsd+ssd}.
Therefore, problem \eqref{expenditureminimizing} is completely solved. 

\begin{remark} Problems \eqref{staticproblem}-\eqref{expenditureminimizing} and their variants for convex, coherent, distortion, or SSD-consistent risk measures have been investigated in Schied \cite{S04}, He and Zhou \cite{HZ11}, F\"ollmer and Schied \cite{FS11}, Mao and Wang \cite{MW20}, Embrechts et al \cite{ESW22}, and Burzoni et al \cite{BMW22}. Risk measure $\Rc$ in this paper lies in none of the aforementioned classes of risk measures. 
\end{remark}

\begin{remark}\label{rmk:emp:fsd}
The classical EMP is as follows: given a quantile function $Q\in\Qsc$, 
\begin{align}\label{opt:emp}
\mini_{X\in L^0}\;\Ebb[\rho X]\quad\mbox{subject to } X\sim Q.
\end{align}
Using the Hardy--Littlewood inequality, it turns out that a
solution to problem \eqref{opt:emp} is given by
$X=Q(1-\xi)$ and the minimum is
$\int_0^1Q(s)Q_{\rho}(1-s)\ds$,
where $\xi\in\Xi$. 
The EMP subject to an FSD constraint or an SSD constraint is as follows:
\begin{align}\label{opt:emp:isd}
\mini_{X\in L^0}\; \Ebb[\rho X]\quad\mbox{subject to } Q_X\in\Qsc_i(Q_0).
\end{align}
For the case of an FSD constraint, i.e., $i=1$,
problem \eqref{opt:emp:isd} is essentially equivalent to problem
\eqref{opt:emp} (with $Q=Q_0$ there) and its solution is thus given by
$X=Q_0(1-\xi)$ and the minimum is
$$x_0=\int_0^1Q_0(s)Q_{\rho}(1-s)\ds,$$
where $\xi\in\Xi$.
For the case of an SSD constraint, i.e., $i=2$,
problem \eqref{opt:emp:isd} has the same minimum $x_0$ as in the case of
an FSD constraint; see, e.g., 
F\"ollmer and Schied \cite[Theorem 3.44]{FS11}. Moreover, a solution of problem \eqref{opt:emp:isd} with $i=2$ is also given by $X=Q_0(1-\xi)$, where $\xi\in\Xi$.
\end{remark}

\section{Proof of Theorem \ref{thm:opt:eqtn}}\label{sec:proof}
To prove Theorem \ref{thm:opt:eqtn}, it suffices to prove the following assertions. 
\begin{enumerate} 
\item $Q^{*}$ is right-continuous.
\item $Q^{*}$ is increasing, which together with the first assertion implies that $Q^{*}$ is a quantile function in $\Qsc$. 
\item $Q^{*}\in\Qsc_1(Q_1)\cap\Qsc_2(Q_2)$. 
\item For any $Q\in\Qsc_1(Q_1)\cap\Qsc_2(Q_2)$, we have $Q\in\Qsc_2(Q^*)$.
\end{enumerate} 
Recall that Assumption \ref{ineq:q12-} guarantees that $f$ is continuous on $[0,1)$. 

\par
Let us start with the first assertion that $Q^{*}$ is right-continuous. For any $t\in(0,1)$, since $\phi\geq f$, there are two possible cases.
\begin{itemize}
\item If $\phi(t)>f(t)$, then by the continuity of $\phi$ and $f$, we have $\phi(t+\Delta)>f(t+\Delta)$ for all sufficiently small $\Delta>0$. Therefore, by definition and the right-continuity of quantiles, 
\[\lim_{\Delta\to 0+}Q^{*}(t+\Delta)=\lim_{\Delta\to 0+}Q_{1}(t+\Delta)=Q_{1}(t)=Q^{*}(t).\]
\item If $\phi(t)=f(t)$, then by definition and the right-continuity of quantiles,
\[\lim_{\Delta\to 0+}Q^{*}(t+\Delta)\leq \lim_{\Delta\to 0+}Q_{1}(t+\Delta)\vee Q_{2}(t+\Delta)=
Q_{1}(t)\vee Q_{2}(t)=Q^{*}(t).\]
By the monotonicity of quantiles, we also have 
\[\lim_{\Delta\to 0+}Q^{*}(t+\Delta) \geq Q^{*}(t).\]
The above two inequalities show that $Q^{*}$ is right-continuous. 
\end{itemize}

{To prove the second and fourth assertions, we first show the following result:}
\begin{align}\label{p2}
\begin{cases} 
\phi'=g,\quad\mbox{a.e. in } (0,1),\\
\phi(0)=0,
\end{cases}
\end{align}
where \[g(t)\trieq (Q_{2}(t)-Q_{1}(t))^{+}\id_{\phi(t)=f(t)}.\] In fact, suppose the ODE in \eqref{opt:equa} holds at $t\in(0,1)$, since $\phi\geq f$, there are two possible cases.
\begin{itemize}
\item If $\phi(t)>f(t)$, then by continuity, we have $\phi(s)>f(s)$ when $s$ is sufficiently close to $t$. It then follows from \eqref{opt:equa} that $\phi$ is a constant near $t$, so $\phi'(t)=0=g(t)$.
\item If $\phi(t)=f(t)$, then by Theorem \ref{thm:ODEsolution}, 
for any $0<\Delta<1-t$, 
\begin{align*} 
\phi(t+\Delta)-\phi(t) &= \max_{ s\leqslant t+\Delta}f(s)-f(t) \geq f (t+\Delta)-f(t), 
\end{align*}
so 
\[\phi'(t)=\liminf_{\Delta\to 0+}\frac{\phi(t+\Delta)-\phi(t)}{\Delta}\geqslant
\liminf_{\Delta\to 0+}\frac{f(t+\Delta)-f(t)}{\Delta}= Q_{2}(t)- Q_{1}(t),\]
by virtue of the right-continuity of quantiles. 
Also, trivially $\phi'\geq 0$, hence $$\phi'(t)\geq \max\{Q_{2}(t)- Q_{1}(t), 0\}=g(t). $$ 
On the other hand, 
\begin{align*} 
\phi(t+\Delta) -\phi(t)
&= \max_{t\leq s\leqslant t+\Delta}f(s)-f(t) \\
&= \max_{ t\leq s\leqslant t+\Delta} \int_{t}^{s} (Q_2(r)-Q_1(r))\dr \\
&\leq \int_{t}^{ t+\Delta} (Q_2(r)-Q_1(r))^{+}\dr, 
\end{align*}
so 
\begin{align*} 
\phi'(t)&=\limsup_{\Delta\to 0+}\frac{\phi(t+\Delta)-\phi(t)}{\Delta}\\
&\leq \limsup_{\Delta\to 0+}\frac{1}{\Delta}\int_{t}^{ t+\Delta} (Q_2(r)-Q_1(r))^{+}\dr \\
&= (Q_2(t)-Q_1(t))^{+}=g(t),
\end{align*} 
thanks to the right-continuity of quantiles. Therefore, $\phi'(t)=g(t)$.
\end{itemize}
Now we have established \eqref{p2}.

\par
{By virtue of \eqref{p2}, } we now prove the second assertion that $Q^{*}$ is increasing. For any $0\leq s<t<1$, as $\phi\geq f$, there are three possible cases.
\begin{itemize}
\item If $\phi(t)=f(t)$, then by definition and the monotonicity of quantiles,
\[Q^*(t)=Q_{1}(t)\vee Q_{2}(t)\geq Q_{1}(s)\vee Q_{2}(s)\geq Q^*(s).\]

\item If $\phi(t)>f(t)$ and $\phi(s)>f(s)$, then by definition and the monotonicity of quantiles,
\[Q^*(t)=Q_{1}(t)\geq Q_{1}(s)=Q^*(s).\]
\item If $\phi(t)>f(t)$ and $\phi(s)=f(s)$, then it follows from \eqref{p2} that 
\begin{align*} 
\int_{s}^{t}Q^*(r)\dr&=\int_{s}^{t}(Q_{1}(r)+g(r))\dr=\phi(t)-\phi(s)+\int_{s}^{t} Q_{1}(r)\dr\\
&>f(t)-f(s)+\int_{s}^{t}Q_{1}(r)\dr=\int_{s}^{t}Q_{2}(r)\dr.\end{align*} 
So there exists some $r\in (s,t)$ such that $Q^*(r)>Q_{2}(r)$. It follows from the definition \eqref{defQ} that $Q_{1}(r)>Q_{2}(r)$, so $Q_{1}(r)= Q_{1}(r)\vee Q_{2}(r)$. Consequently, 
\[Q^*(t)=Q_{1}(t)\geq Q_{1}(r)= Q_{1}(r)\vee Q_{2}(r)\geq Q_{1}(s)\vee Q_{2}(s)\geq Q^*(s)\]
by the monotonicity of quantiles. 
\end{itemize}
We have now proved the second assertion. 
\par 
{Thanks to \eqref{p2} again, } we can now show the third assertion that $Q^{*}\in\Qsc_1(Q_1)\cap\Qsc_2(Q_2)$. By \eqref{opt:equa} and \eqref{p2}, we see that for a.e. $t\in(0,1)$, 
\begin{equation*} 
\min\left\{Q^{*}(t)-Q_1(t), \int_0^tQ^{*}(r)\dr- \int_0^tQ_{2}(r)\dr\right\} 
=\min\left\{g(t), \int_0^tg(r)\dr-f(t)\right\}=0.
\end{equation*}
By the right-continuity of quantiles and integrals, the above holds for all $t\in(0,1)$. Hence, $Q^{*}\in\Qsc_1(Q_1)\cap\Qsc_2(Q_2)$. 
\par
We now show the fourth assertion that $Q\in\Qsc_2(Q^*)$ for any $Q\in\Qsc_1(Q_1)\cap\Qsc_2(Q_2)$. In fact, if $Q\in\Qsc_1(Q_1)\cap\Qsc_2(Q_2)$, then 
\begin{equation} \label{equQ2}
\min\left\{Q(t)-Q_1(t), \int_0^tQ(r)\dr-\int_0^tQ_{2}(r)\dr\right\}\ge0,\quad\mbox{for all } t\in(0,1).
\end{equation}
The following comparison principle shows that $\int_0^tQ(r)\dr\geq \int_0^tQ^{*}(r)\dr$ for all $t\in(0,1)$. Hence, 
$Q\in\Qsc_2(Q^*)$, which completes the proof of Theorem \ref{thm:opt:eqtn}.

\begin{lemma}[Comparison principle]\label{comparison}
Assume $\varphi_{1}$ and $\varphi_{2}$ are absolutely continuous functions that satisfy the following variational inequalities on $[0,1)$, respectively:
\begin{align*}
\min\{\varphi'_{1}-f_{1},\varphi_{1}-g_{1}\}\geq0 \mbox{ a.e.},
\quad
\min\{\varphi'_{2}-f_{2},\varphi_{2}-g_{2}\}\leq0\mbox{ a.e.}
\end{align*}
If $f_{1}\geq f_{2}$, $g_{1}\geq g_{2}$ a.e. on $[0,1)$, and $\varphi_{1}(0)\geq \varphi_{2}(0)$, then $\varphi_{1}\geq \varphi_{2}$ on $[0,1)$.
\end{lemma}
\proof
Suppose, on the contrary, there exists a $t\in(0,1)$ such that $\varphi_{1}(t)<\varphi_{2}(t)$. Let 
\[s=\sup\{0\leq r<t\mid \varphi_{1}(r)\geq \varphi_{2}(r)\}.\] Then by the continuity of $\varphi_{1}$ and $\varphi_{2}$ and $\varphi_{1}(0)\geq \varphi_{2}(0)$, we get $0\leq s<t$, $\varphi_{1}(s)\geq\varphi_{2}(s)$ and $\varphi_{1}(r)<\varphi_{2}(r)$ for all $r\in(s,t]$. This together with the first variational inequality and the fact $g_{1}\geq g_{2}$ gives 
\[\varphi_{2}>\varphi_{1}\geq g_{1}\geq g_{2} \quad \mbox{a.e. in } (s,t].\]
Then, by the two variational inequalities and the fact $f_{1}\geq f_{2}$, we get
\[\varphi'_{2}\leq f_{2}\leq f_{1}\leq \varphi'_{1}\quad\mbox{a.e. in } (s,t]. \]
It follows 
\[\varphi_{2}(t)-\varphi_{2}(s)=\int_{s}^{t}\varphi'_{2}(r)\dr\leq \int_{s}^{t}\varphi'_{1}(r)\dr=\varphi_{1}(t)-\varphi_{1}(s),\]
contradicting the fact that $\varphi_{1}(t)<\varphi_{2}(t)$ and $\varphi_{1}(s)\geq\varphi_{2}(s)$.
\qed

\bibliographystyle{siamplain}
\bibliography{references}

\begin{thebibliography}{}


\bibitem{BBM20} Bignozzi, V., Burzoni, M., and Munari, C. (2020): Risk measures based on benchmark loss distributions, \textit{Journal of Risk and Insurance} \textbf{87}(2): 437--475.

\bibitem{BMW22} Burzoni, M., C. Munari and R. Wang (2022): Adjusted Expected Shortfall, \textit{Journal of Banking and Finance}, forthcoming.

\bibitem{CD06} Carlier, G. and R.-A. Dana (2006): 
Law Invariant Concave Utility Functions and Optimization Problems with Monotonicity and
Comonotonicity Constraints, \textit{Stat. Decis.} \textbf{24},
127--152.


\bibitem{D88} Dybvig, P. H. (1988): Distributional Analysis of Portfolio Choice, \textit{J. Business}
\textbf{61}, 369--398.

\bibitem{ESW22} Embrechts, P., A. Schied and R. Wang (2022): Robustness in the optimization of risk measures, 
\textit{Operations Research} \textbf{70}, 95--110. 
\bibitem{FS11} F\"ollmer, H. and A. Schied (2011): \textit{Stochastic
Finance: An Introduction in Discrete Time (3rd Edition)}. Berlin:
Walter de Gruyter.

\bibitem{HZ11} He, X. D. and X. Y. Zhou (2011): Portfolio Choice via
Quantiles, \textit{Math. Finance} \textbf{21}, 203--231.

\bibitem{JZ08} Jin, H. and X. Y. Zhou (2008): Behavioral Portfolio
Selection in Continuous Time, \textit{Math. Finance} \textbf{18},
385--426. 

\bibitem{LM22} Liebrich, F.-B., and C. Munari (2022): Law-invariant functionals that collapse to the Mean: beyond convexity, {\it Mathematics and Financial Economics} {\bf 16}, 447--480. 

\bibitem{MW20} Mao, T., and Wang, R. (2020): Risk aversion in regulatory capital calculation, {\it SIAM Journal on Financial Mathematics} {\bf 11}, 169--200.

\bibitem{MR01} M\"uller, A. and R\"uschendorf L. (2001): On the optimal stopping values induced by general depence structures, {\it J. Appl. Probab.} {\bf 38}, 672--684. 

\bibitem{RY99} Revuz, D. and M. Yor (1999): \textit{Continuous Martingales and Brownian Motion (3rd Edition)}. New York: Springer.

\bibitem{S04} Schied, A (2004): On the Neyman-Pearson Problem for Law-Invariant Risk Measures and
Robust Utility Functionals, \textit{Ann. Appl. Probab.}
\textbf{14}, 1398--1423.

\bibitem{WX21} Wang, X. and J. Xia (2021): Expected Utility Maximization with Stochastic Dominance Constraints in Complete Markets, \textit{SIAM Journal on Financial Mathematics} \textbf{12}, 1054--1111.

\bibitem{XZ16} Xia, J. and X. Y. Zhou (2016): Arrow-Debreu Equilibria for Rank-Dependent
Utilities, \textit{Math. Finance} \textbf{26}, 558--588.

\bibitem{X14} Xu, Z. Q. (2014): A New Characterization of Comonotonicity and Its Application in Behavioral Finance, \textit{J. Math. Anal. Appl.} \textbf{418}, 612--625. 

\bibitem{X16} Xu, Z. Q. (2016): A Note on the Quantile Formulation, \textit{Math. Finance} \textbf{26}, 589--601.



\end{thebibliography}

\end{document}